\documentclass{amsart}

\theoremstyle{definition}

\theoremstyle{remark}

\numberwithin{equation}{section}

\begin{document}

\title[Double difference property for H\"{o}lder continuous functions]{The double difference property for the class of locally H\"{o}lder continuous functions}

\author{Rashid A. Aliev}

\address{Institute of Mathematics and Mechanics, NAS of Azerbaijan, Baku, Azerbaijan \newline Baku State University, Baku, Azerbaijan}

\email{aliyevrashid@mail.ru}

\author{Aysel A. Asgarova}

\address{Azerbaijan University of Languages, Baku, Azerbaijan}

\email{asgarova2016@mail.ru}

\author{Vugar E. Ismailov}

\address{Institute of Mathematics and Mechanics, NAS of Azerbaijan, Baku, Azerbaijan}

\email{vugaris@mail.ru}



\begin{abstract}
In this paper, we show that the pair of classes of locally H\"{o}lder continuous functions
(considered on $\mathbb{R}$ and $\mathbb{R}^{2}$, respectively) has the double difference property.
\end{abstract}

\maketitle

\section{Introduction}

The notions \textit{difference property} and \textit{double difference property}
are due to de Bruijn \cite{1} and Laczkovich \cite{5},
respectively. These properties for various classes of real functions were
investigated by many authors. We refer the reader to Laczkovich's survey
paper \cite{6} for a detailed source of information on this topic.

For a fixed function $g:\mathbb{R}\rightarrow \mathbb{R}$ and any $h\in
\mathbb{R}$ we define the difference function $\Delta _{h}g:\mathbb{R}%
\rightarrow \mathbb{R}$ by
\begin{equation*}
\Delta _{h}g(x)=g(x+h)-g(x)
\end{equation*}%
and the double difference function $Dg:$ $\mathbb{R}^{2}\rightarrow \mathbb{R%
}$ by
\begin{equation*}
Dg(x,y)=g(x+y)-g(x)-g(y).
\end{equation*}

Let $\mathcal{F}$ be a class of functions defined on $\mathbb{R}$ and $%
\mathcal{F}_{2}$ be a class of functions defined on $\mathbb{R}^{2}$. The
class $\mathcal{F}$ is said to have the difference property if every
function $g:\mathbb{R}\rightarrow \mathbb{R}$, for which $\Delta _{h}g\in
\mathcal{F}$ for each $h\in \mathbb{R}$, is of the form $g=f+A$, where $f\in
\mathcal{F}$ and $A$ is an additive function (see \cite{1}). A function $A$
is called additive if it satisfies the Cauchy functional equation $%
A(x+y)=A(x)+A(y)$. The pair $\left( \mathcal{F},\mathcal{F}_{2}\right) $ is
said to have the double difference property if whenever $Dg\in \mathcal{F}%
_{2}$ holds for a function $g:\mathbb{R}\rightarrow \mathbb{R}$, then $g$ is
of the form $g=f+A$, where $f\in \mathcal{F}$ and $A$ is additive (see \cite%
{5}).

de Bruijn \cite{1} was the first who showed that the class of continuous
functions has the difference property and thus resolved Erd\"{o}s's famous
conjecture. He also proved that the difference property holds for a large
number of essential function classes (see \cite{1,2}). Some of these classes
are

\smallskip

1) $C^{k}(\mathbb{R)}$, functions with continuous derivatives up to order $k$%
;

2) $C^{\infty }(\mathbb{R)}$, infinitely differentiable functions;

3) analytic functions;

4) functions which are absolutely continuous on any finite interval;

5) functions having bounded variation over any finite interval;

6) algebraic polynomials;

7) trigonometric polynomials;

8) Riemann integrable functions.

\smallskip

However, the class $\mathcal{L}$ of Lebesgue measurable functions fails to
have this property if we assume the continuum hypothesis (see \cite{1,5}).
It was conjectured by Erd\"{o}s that every function $g:\mathbb{R}\rightarrow
\mathbb{R}$ for which $\Delta _{h}g(x)$ is measurable for each $h$, is of
the form $g=f+A+S$, where $f$ is measurable, $A$ is additive and $S$ has the
property that $\Delta _{h}S(x)=0$ for almost all $x$. Laczkovich \cite{5}
solved this conjecture affirmatively and moreover proved that the pair $(%
\mathcal{L},\mathcal{L}_{2})$ has the double difference property, where $%
\mathcal{L}_{2}$ denotes the class of Lebesgue measurable functions defined
on $\mathbb{R}^{2}$. It was also proved in \cite{5} that the double
difference property holds for Baire $\alpha $ functions. Later Tabor and
Tabor \cite{7} proved that the class $C^{n}(X,Y)$ of $n$-times continuously
differentiable functions defined on a real normed space $X$ and taking
values in a real Banach space $Y$ has the double difference property.
Kotlicka \cite{3} showed that several pairs of classes of functions have the
double difference property. Among them there are approximately continuous
functions, pointwise continuous functions, essentially continuous functions
(considered on $\mathbb{R}$ and $\mathbb{R}^{2}$, respectively) and $L_{p}$%
-classes for $0<p<\infty $ (considered on $\mathbb{T}$ and $\mathbb{T}^{2}$,
respectively, where $\mathbb{T}$ is the torus).

In \cite{T1}, Tabor proved that the pair of classes of Lipschitz functions
defined on a metric semigroup $G$ and $G\times G$, respectively, with values
in a reflexive Banach space $E$ has the double difference property. For
finite dimensional Banach spaces $X$ and $Y$, Tabor and Tabor \cite{T2}
showed that the double difference property holds for the pair of classes of $%
Y$-valued Lipschitz functions defined on a convex set $K\subset X$ such that
$0\in K$ $\ $and on the set
\begin{equation*}
C(K)=\left\{ (x,y)\in X\times X:x\in K,y\in K,x+y\in K\right\} ,
\end{equation*}%
respectively. Consequently, the double difference property holds for the
pair of classes of real Lipschitz functions defined on an interval $I$
containing zero and on the set $C(I)$, respectively.

In this paper, we prove that for any $\alpha \in (0,1]$ the pair of classes
of locally H\"{o}lder real continuous functions (considered on $\mathbb{R}$
and $\mathbb{R}^{2}$, respectively) with exponent $\alpha $ has the double
difference property.

\bigskip

\section{Main result}

We start this section with the definition of modulus of continuity of a
multivariate function and some notation. Let $f(\mathbf{x}%
)=f(x_{1},...,x_{s})$, $s\geq 1,$ be any $s$-variable function defined on a
set $\Omega \subset \mathbb{R}^{s}$. The function
\begin{equation*}
\omega (f;\delta ;\Omega )=\sup \left\{ \left\vert f(\mathbf{x})-f(\mathbf{y}%
)\right\vert :\mathbf{x},\mathbf{y}\in \Omega ,\text{ }\left\vert \mathbf{x}-%
\mathbf{y}\right\vert \leq \delta \right\} ,\text{ }0\leq \delta \leq
diam\Omega ,
\end{equation*}%
is called the modulus of continuity of $f$ on $\Omega $. We will also use
the notation $\omega _{\mathbb{Q}}(f;\delta ;\Omega )$, which stands for the
function $\omega (f;\delta ;\Omega \cap \mathbb{Q}^{s})$. Here $\mathbb{Q}%
^{s}$ denotes the space of $s$-dimensional vectors with rationale
coordinates. Clearly, $\omega _{\mathbb{Q}}(f;\delta ;\Omega )$ makes sense
if the set $\Omega \cap \mathbb{Q}^{s}$ is not empty. Note that we always
have the inequality $\omega _{\mathbb{Q}}(f;\delta ;\Omega )\leq \omega
(f;\delta ;\Omega )$ and the strong equality $\omega _{\mathbb{Q}}(f;\delta
;\Omega )=\omega (f;\delta ;\Omega )$ holds for continuous $f$ and certain
sets $\Omega $. For example, this holds if for any $\mathbf{x},\mathbf{y}\in
\Omega $ with $\left\vert \mathbf{x}-\mathbf{y}\right\vert \leq \delta $
there exist sequences $\left\{ \mathbf{x}_{n}\right\} ,\left\{ \mathbf{y}%
_{n}\right\} \subset \Omega \cap \mathbb{Q}^{s}$ such that $\mathbf{x}%
_{n}\rightarrow \mathbf{x}$, $\mathbf{y}_{n}\rightarrow \mathbf{y}$ and $%
\left\vert \mathbf{x}_{n}-\mathbf{y}_{n}\right\vert \leq \delta ,$ for all $%
n $. There are many sets $\Omega $, which satisfy this property.

The class $H_{\alpha }^{\left( loc\right) }\left( \mathbb{R}^{s}\right) $ of
locally H\"{o}lder continuous functions with exponent $\alpha $ is defined
as the class of functions $f$ for which $\omega (f;\delta ;\Omega )\leq
K\delta ^{\alpha }$ for any compact set $\Omega \subset \mathbb{R}^{s}$.
Here $K$ depends on $\Omega $.

\bigskip

Our main result is the following theorem.

\bigskip

\textbf{Theorem 2.1.} \textit{Assume a function $g:\mathbb{R}\rightarrow
\mathbb{R}$ is such that the bivariate function $\mathit{g(x+y)}-g(x)-g(y)$
is locally H\"{o}lder continuous with exponent $\alpha$. Then there exist a
function $f\in H_{\alpha }^{\left( loc\right) }\left( \mathbb{R}\right) $
and an additive function $A$ such that $g=f+A$.}

\bigskip

To prove this theorem we need the following auxiliary lemma.

\bigskip

\textbf{Lemma 2.1.} \textit{Assume a function $F\in C(\mathbb{R}^{2})$ has
the form}
\begin{equation}
F(x,y)=g(x+y)-g(x)-g(y),  \label{3.1}
\end{equation}%
\textit{where $g$ is an arbitrarily behaved function. Then the following
inequality holds}

\begin{equation}
\omega _{\mathbb{Q}}(g;\delta ;[-M,M])\leq 2\delta \left\vert
g(1)-g(0)\right\vert +3\omega \left( F;\delta ;[-M,M]^{2}\right) ,
\label{3.2}
\end{equation}%
\textit{where $\delta \in \left( 0,\frac{1}{2}\right) \cap \mathbb{Q}$ and $%
M\geq 1$.}

\bigskip

\begin{proof} Consider the function $h(t)=g(t)-g(0)$ and write \eqref{3.1}
in the form

\begin{equation}
G(x,y)=h(x+y)-h(x)-h(y),  \label{3.3}
\end{equation}%
where

\begin{equation*}
G(x,y)=F(x,y)+g(0).
\end{equation*}%
Note that the functions $g$ and $h$, as well as the functions $F$ and $G,$
have the common modulus of continuity. Thus we prove the lemma if we prove
it for the pair $\left\langle G,h\right\rangle .$

\bigskip Since $h(0)=0,$ it follows from \eqref{3.3} that

\begin{equation}
G(x,0)=G(0,y)=0.  \label{3.4}
\end{equation}%
Obviously, for any real number $x,$
\begin{eqnarray*}
G(x,x) &=&h(2x)-2h(x); \\
G(x,2x) &=&h(3x)-h(x)-h(2x); \\
&&\cdot \cdot \cdot \\
G(x,(k-1)x) &=&h(kx)-h(x)-h((k-1)x).
\end{eqnarray*}

We obtain from the above equalities that

\begin{eqnarray*}
h(2x) &=&2h(x)+G(x,x), \\
h(3x) &=&3h(x)+G(x,x)+G(x,2x), \\
&&\cdot \cdot \cdot \\
h(kx) &=&kh(x)+G(x,x)+G(x,2x)+\cdot \cdot \cdot +G(x,(k-1)x).
\end{eqnarray*}%
Thus for any nonnegative integer $k$,

\begin{equation}
h(x)=\frac{1}{k}h(kx)-\frac{1}{k}\left[ G(x,x)+G(x,2x)+\cdot \cdot \cdot
+G(x,(k-1)x)\right] .  \label{3.5}
\end{equation}

Consider now the simple fraction $\frac{p}{n}\in (0,\frac{1}{2})$ and set $%
m_{0}=\left[ \frac{n}{p}\right] .$ Here $[r]$ denotes the whole number part
of $r$. Clearly, $m_{0}\geq 2$ and the remainder $p_{1}=n-m_{0}p<p.$ Taking $%
x=\frac{p}{n}$ and $k=m_{0}$ in \eqref{3.5} gives us the following equality

\begin{equation*}
h\left( \frac{p}{n}\right) =\frac{1}{m_{0}}h\left( 1-\frac{p_{1}}{n}\right)
\end{equation*}%
\begin{equation}
-\frac{1}{m_{0}}\left[ G\left( \frac{p}{n},\frac{p}{n}\right) +G\left( \frac{%
p}{n},\frac{2p}{n}\right) +\cdot \cdot \cdot +G\left( \frac{p}{n},(m_{0}-1)%
\frac{p}{n}\right) \right] .  \label{3.6}
\end{equation}%
On the other hand, since

\begin{equation*}
G\left( \frac{p_{1}}{n},1-\frac{p_{1}}{n}\right) =h(1)-h\left( \frac{p_{1}}{n%
}\right) -h\left( 1-\frac{p_{1}}{n}\right) ,
\end{equation*}%
it follows from \eqref{3.6} that

\begin{equation*}
h\left( \frac{p}{n}\right) =\frac{h(1)}{m_{0}}-\frac{1}{m_{0}}\left[ G\left(
\frac{p}{n},\frac{p}{n}\right) +\cdot \cdot \cdot +G\left( \frac{p}{n}%
,(m_{0}-1)\frac{p}{n}\right) +G\left( \frac{p_{1}}{n},1-\frac{p_{1}}{n}%
\right) \right]
\end{equation*}%
\begin{equation}
-\frac{1}{m_{0}}h\left( \frac{p_{1}}{n}\right) .  \label{3.7}
\end{equation}

Put $m_{1}=\left[ \frac{n}{p_{1}}\right] $, $p_{2}=n-m_{1}p_{1}.$ Clearly, $%
0\leq p_{2}<p_{1}$. Similar to \eqref{3.7}, we can write that

\begin{equation*}
h\left( \frac{p_{1}}{n}\right) =\frac{h(1)}{m_{1}}-\frac{1}{m_{1}}\left[
G\left( \frac{p_{1}}{n},\frac{p_{1}}{n}\right) +\cdot \cdot \cdot +G\left(
\frac{p_{1}}{n},(m_{1}-1)\frac{p_{1}}{n}\right) +G\left( \frac{p_{2}}{n},1-%
\frac{p_{2}}{n}\right) \right]
\end{equation*}%
\begin{equation}
-\frac{1}{m_{1}}h\left( \frac{p_{2}}{n}\right) .  \label{3.8}
\end{equation}

Let us make a convention that \eqref{3.7} is the $1$-st and \eqref{3.8} is
the $2$-nd formula. One can continue this process by defining the chain of
pairs $(m_{2},p_{3}),$ $(m_{3},p_{4})$ until the pair $(m_{k-1},p_{k})$ with
$p_{k}=0$ and writing out the corresponding formulas for each pair. For
example, the last $k$-th formula will be of the form
\begin{equation*}
h\left( \frac{p_{k-1}}{n}\right) =\frac{h(1)}{m_{k-1}}
\end{equation*}%
\begin{equation*}
-\frac{1}{m_{k-1}}\left[ G\left( \frac{p_{k-1}}{n},\frac{p_{k-1}}{n}\right)
+\cdot \cdot \cdot +G\left( \frac{p_{k-1}}{n},(m_{k-1}-1)\frac{p_{k-1}}{n}%
\right) +G\left( \frac{p_{k}}{n},1-\frac{p_{k}}{n}\right) \right]
\end{equation*}%
\begin{equation}
-\frac{1}{m_{k-1}}h\left( \frac{p_{k}}{n}\right) .  \label{3.9}
\end{equation}%
Note that in \eqref{3.9}, $h\left( \frac{p_{k}}{n}\right) =0$ and $G\left(
\frac{p_{k}}{n},1-\frac{p_{k}}{n}\right) =0$. Considering now the $k$-th
formula in the $(k-1)$-th formula, then the obtained formula in the $(k-2)$%
-th formula, and so forth, we will finally arrive at the equality

\begin{equation*}
h\left( \frac{p}{n}\right) =h(1)\left[ \frac{1}{m_{0}}-\frac{1}{m_{0}m_{1}}%
+\cdot \cdot \cdot +\frac{(-1)^{k-1}}{m_{0}m_{1}\cdot \cdot \cdot m_{k-1}}%
\right]
\end{equation*}

\begin{equation*}
-\frac{1}{m_{0}}\left[ G\left( \frac{p}{n},\frac{p}{n}\right) +\cdot \cdot
\cdot +G\left( \frac{p}{n},(m_{0}-1)\frac{p}{n}\right) +G\left( \frac{p_{1}}{%
n},1-\frac{p_{1}}{n}\right) \right]
\end{equation*}

\begin{equation*}
+\frac{1}{m_{0}m_{1}}\left[ G\left( \frac{p_{1}}{n},\frac{p_{1}}{n}\right)
+\cdot \cdot \cdot +G\left( \frac{p_{1}}{n},(m_{1}-1)\frac{p_{1}}{n}\right)
+G\left( \frac{p_{2}}{n},1-\frac{p_{2}}{n}\right) \right]
\end{equation*}

\begin{equation*}
+\cdot \cdot \cdot +
\end{equation*}%
\begin{equation}
\frac{(-1)^{k}}{m_{0}m_{1}\cdot \cdot \cdot m_{k-1}}\left[ G\left( \frac{%
p_{k-1}}{n},\frac{p_{k-1}}{n}\right) +\cdot \cdot \cdot +G\left( \frac{%
p_{k-1}}{n},(m_{k-1}-1)\frac{p_{k-1}}{n}\right) \right] .  \label{3.10}
\end{equation}%
Taking into account \eqref{3.4}, we obtain from \eqref{3.10} that

\begin{equation*}
\left\vert h\left( \frac{p}{n}\right) \right\vert \leq \left[ \frac{1}{m_{0}}%
-\frac{1}{m_{0}m_{1}}+\cdot \cdot \cdot +\frac{(-1)^{k-1}}{m_{0}m_{1}\cdot
\cdot \cdot m_{k-1}}\right] \left\vert h(1)\right\vert
\end{equation*}

\begin{equation}
+\left[ 1+\frac{1}{m_{0}}+\cdot \cdot \cdot +\frac{1}{m_{0}\cdot \cdot \cdot
m_{k-2}}\right] \omega \left( G;\frac{p}{n};[0,1]^{2}\right).  \label{3.11}
\end{equation}%
Since $m_{0}\leq m_{1}\leq \cdot \cdot \cdot \leq m_{k-1},$ it is not
difficult to see that in \eqref{3.11}

\begin{equation*}
\frac{1}{m_{0}}-\frac{1}{m_{0}m_{1}}+\cdot \cdot \cdot +\frac{(-1)^{k-1}}{%
m_{0}m_{1}\cdot \cdot \cdot m_{k-1}}\leq \frac{1}{m_{0}}
\end{equation*}%
and

\begin{equation*}
1+\frac{1}{m_{0}}+\cdot \cdot \cdot +\frac{1}{m_{0}\cdot \cdot \cdot m_{k-2}}%
\leq \frac{m_{0}}{m_{0}-1}.
\end{equation*}%
Considering the above two inequalities in \eqref{3.11} we obtain that

\begin{equation}
\left\vert h\left( \frac{p}{n}\right) \right\vert \leq \frac{\left\vert
h(1)\right\vert }{m_{0}}+\frac{m_{0}}{m_{0}-1}\omega \left( G;\frac{p}{n}%
;[0,1]^{2}\right) .  \label{3.12}
\end{equation}%
Since $m_{0}=\left[ \frac{n}{p}\right] \geq 2,$ it follows from \eqref{3.12}
that

\begin{equation}
\left\vert h\left( \frac{p}{n}\right) \right\vert \leq \frac{2p\left\vert
h(1)\right\vert }{n}+2\omega \left( G;\frac{p}{n};[0,1]^{2}\right) .
\label{3.13}
\end{equation}

Let now $\delta \in \left( 0,\frac{1}{2}\right) \cap \mathbb{Q}$ be a
rational increment, $M\geq 1$ and $x,x+\delta $ be two points in $\left[ -M,M%
\right] \cap \mathbb{Q}.$ By \eqref{3.3}, \eqref{3.4} and \eqref{3.13} we
can write that

\begin{equation}
\left\vert h(x+\delta )-h(x)\right\vert \leq \left\vert h(\delta
)\right\vert +\left\vert G(x,\delta )\right\vert \leq 2\delta \left\vert
h(1)\right\vert +3\omega \left( G;\delta ;[-M,M]^{2}\right) .  \label{3.14}
\end{equation}%
Now \eqref{3.2} follows from \eqref{3.14} and the definitions of $h$ and $G$.
\end{proof}

\bigskip

\textit{Remark 1.} Under the assumptions of Lemma 2.1, the restriction of $g$
to the set of rational numbers is uniformly continuous on any interval $%
[-M,M]\cap \mathbb{Q}$ and hence continuous on $\mathbb{Q}$.

\bigskip

Now we are ready to prove Theorem 2.1.

\begin{proof} Let us put

\begin{equation}
F(x,y)=g(x+y)-g(x)-g(y)  \label{3.00}
\end{equation}%
and consider the function

\begin{equation*}
u(t)=g(t)-\left[ g(1)-g(0)\right] t.
\end{equation*}%
Obviously, $u(1)=u(0)$ and

\begin{equation}
F(x,y)=u(x+y)-u(x)-u(y).  \label{3.17}
\end{equation}%
By Lemma 2.1, the restriction of $u$ to $\mathbb{Q}$ is continuous and
uniformly continuous on every interval $[-M,M]\cap \mathbb{Q}$. Denote this
restriction by $v$.

Let $y$ be any real number and $\{y_{n}\}_{n=1}^{\infty }$ be any sequence
of rational numbers converging to $y$. We can choose $M>0$ so that $y_{n}\in
\lbrack -M,M]$ for any $n\in \mathbb{N}$. It follows from the uniform
continuity of $v$ on $[-M,M]\cap \mathbb{Q}$ that the sequence $%
\{v(y_{n})\}_{n=1}^{\infty }$ is Cauchy. Thus there exits a finite limit $%
\lim_{n\rightarrow \infty }v(y_{n})$. It is not difficult to see that this
limit does not depend on the choice of $\{y_{n}\}_{n=1}^{\infty }$.

Let $f$ denote the following extension of $v$ to the set of real numbers.

\begin{equation*}
f(y)=\left\{
\begin{array}{c}
v(y),\text{ if }y\in \mathbb{Q}\text{;} \\
\lim_{n\rightarrow \infty }v(y_{n}),\text{ if }y\in \mathbb{R}\backslash
\mathbb{Q}\text{ and }\{y_{n}\}\text{ is a sequence in }\mathbb{Q}\text{
tending to }y.%
\end{array}%
\right.
\end{equation*}%
In view of the above arguments, $f$ is well defined on the whole real line.
Let us prove that $f$ is the function we seek.

Consider an arbitrary point $(x,y)\in \mathbb{R}^{2}$ and a sequence of
points $\{(x_{n},y_{n})\}_{n=1}^{\infty }$ with rationale coordinates
tending to $(x,y)$. Taking into account \eqref{3.17}, we can write that

\begin{equation}
F(x_{n},y_{n})=v(x_{n}+y_{n})-v(x_{n})-v(y_{n}),\text{ for all }n=1,2,...,
\label{3.18}
\end{equation}%
since $v$ is the restriction of $u$ to $\mathbb{Q}$. Tending $n\rightarrow
\infty $ in both sides of \eqref{3.18} we obtain that

\begin{equation}
F(x,y)=f(x+y)-f(x)-f(y).  \label{3.0}
\end{equation}

Set $A=g-f$. It follows from \eqref{3.00} and \eqref{3.0} that $A$ is
additive. Let us now prove that $f\in H_{\alpha }^{\left( loc\right) }\left(
\mathbb{R}\right) $. Since $v(1)=v(0)$ we obtain from \eqref{3.17} and %
\eqref{3.2} that for $\delta \in \left( 0,\frac{1}{2}\right) \cap \mathbb{Q}$%
, $M\geq 1$ and any numbers $a,b\in \lbrack -M,M]\cap \mathbb{Q}$, $%
\left\vert a-b\right\vert \leq \delta ,$ the following inequality holds

\begin{equation}
\left\vert v(a)-v(b)\right\vert \leq 3\omega \left( F;\delta
;[-M,M]^{2}\right) .  \label{3.19}
\end{equation}%
Consider now any real numbers $r_{1}$ and $r_{2}$ satisfying $r_{1},r_{2}\in
\lbrack -M,M]$, $\left\vert r_{1}-r_{2}\right\vert \leq \delta $ and take
sequences $\{a_{n}\}_{n=1}^{\infty }\subset \lbrack -M,M]\cap \mathbb{Q}$, $%
\{b_{n}\}_{n=1}^{\infty }\subset \lbrack -M,M]\cap \mathbb{Q}$ with the
property $\left\vert a_{n}-b_{n}\right\vert \leq \delta ,$ $n=1,2,...,$ and
tending to $r_{1}$ and $r_{2}$, respectively. By \eqref{3.19},
\begin{equation*}
\left\vert v(a_{n})-v(b_{n})\right\vert \leq 3\omega \left( F;\delta
;[-M,M]^{2}\right) .
\end{equation*}%
If we take limits\ on both sides of the above inequality, we obtain that

\begin{equation}
\left\vert f(r_{1})-f(r_{2})\right\vert \leq 3\omega \left( F;\delta
;[-M,M]^{2}\right) ,  \label{3.20}
\end{equation}%
which means that $f$ is uniformly continuous on $[-M,M]$ and hence it is
continuous on the whole real line.

It follows from \eqref{3.20} that

\begin{equation}
\omega \left( f;\delta ;[-M,M]\right) \leq 3\omega \left( F;\delta
;[-M,M]^{2}\right) .  \label{3.21}
\end{equation}%
Note that in \eqref{3.21} $\delta $ is a rational number from the interval $%
\left( 0,\frac{1}{2}\right) $. Since the modulus of continuity of a
continuous function is continuous from the right (see \cite{Kol}), it
follows that, in fact, \eqref{3.21} is valid for all $\delta \in \left[ 0,%
\frac{1}{2}\right) $. Since $F\in H_{\alpha }^{\left( loc\right) }\left(
\mathbb{R}^{2}\right) $, we obtain from \eqref{3.21} that

\begin{equation}
\omega \left( f;\delta ;[-M,M]\right) \leq K\delta ^{\alpha },\text{ where }%
0\leq \delta <\frac{1}{2}\text{.}  \label{3.22}
\end{equation}

Let now $\frac{1}{2}\leq \delta \leq 2M$. We can write that

\begin{equation}
\omega \left( f;\delta ;[-M,M]\right) \leq 2\left\Vert f\right\Vert
_{C([-M,M])}\leq 2^{1+\alpha }\left\Vert f\right\Vert _{C([-M,M])}\delta
^{\alpha },\text{ where }\frac{1}{2}\leq \delta \leq 2M\text{.}  \label{3.23}
\end{equation}

The inequalities \eqref{3.22} and \eqref{3.23} show that $f$ is H\"{o}lder
continuous on $[-M,M]$ with exponent $\alpha $. Since $M$ is an arbitrary
number not less than $1$ and any compact $\Omega $ is contained in a closed
interval of the form $[-M,M]$, we obtain that $f\in H_{\alpha }^{\left(
loc\right) }\left( \mathbb{R}\right) $.
\end{proof}

\textit{Remark 2.} The above proof shows that for any compact set $\Omega
\subset $ $\mathbb{R}$ the pair of H\"{o}lder continuous function classes $%
H_{\alpha }\left( \Omega \right) $ and $H_{\alpha }\left( \Omega \times
\Omega \right) $ has the double difference property. This holds, in
particular, for the pair of classes of Lipschitz functions defined on $%
\Omega $ and $\Omega \times \Omega $, respectively. The last assertion
complements the corresponding result of J. Tabor and J. Tabor \cite{T2} in
the real space setting (see Introduction).

\bigskip

\textit{Remark 3.} Theorem 2.1 is not only an existence result. It's proof
gives a recipe for constructing the function $f$. It also allows us to
estimate the modulus of continuity of $f$ in terms of the modulus of
continuity of $g(x+y)-g(x)-g(y)$ (see (2.21)).

\end{document}